# Adaptive decision-making for stochastic service network design

**Javier Durán-Micco[a*], Bilge Atasoy[a]**

[a]Maritime and Transport Technology (MTT), Mechanical Engineering (ME), Delft University of Technology/ Mekelweg 2, 2628CD, Delft, The Netherlands.
*Corresponding author. Email: j.duranmicco@tudelft.nl.

## Abstract

This paper addresses the Service Network Design (SND) problem for a logistics service provider (LSP) operating in a multimodal freight transport network, considering uncertain travel times and limited truck fleet availability. A two-stage optimization approach is proposed, which combines metaheuristics, simulation and machine learning components. This solution framework integrates tactical decisions, such as transport request acceptance and capacity booking for scheduled services, with operational decisions, including dynamic truck allocation, routing, and re-planning in response to disruptions. A simulated annealing (SA) metaheuristic is employed to solve the tactical problem, supported by an adaptive surrogate model trained using a discrete-event simulation model that captures operational complexities and cascading effects of uncertain travel times. The performance of the proposed method is evaluated using benchmark instances. First, the SA is tested on a deterministic version of the problem and compared to state-of-the-art results, demonstrating it can improve the solution quality and significantly reduce the computational time. Then, the proposed SA is applied to the more complex stochastic problem. Compared to a benchmark algorithm that executes a full simulation for each solution evaluation, the learning-based SA generates high quality solutions while significantly reducing computational effort, achieving only a 5% difference in objective function value while cutting computation time by up to 20 times. These results demonstrate the strong performance of the proposed algorithm in solving complex versions of the SND. Moreover, they highlight the effectiveness of integrating diverse modeling and optimization techniques, and the potential of such approaches to efficiently address freight transport planning challenges.

# 1 Introduction

Increasing demand for more efficient, reliable, and sustainable freight transport systems have multiplied the need for advanced planning, modelling, and optimization techniques. Logistics Service Providers (LSPs) play a key role in such systems, acting as intermediaries who organize the transport of goods among various actors. Through the consolidation of cargo from multiple shippers and utilization of the services of multiple carriers, LSPs can achieve economies of scale and improved resource utilization. This central coordination role, however, introduces significant planning complexity, particularly because LSPs must make decisions on a multitude of levels and under many sources of uncertainty.

One of the most significant tactical problems LSPs must deal with is the Service Network Design (SND) problem, which addresses the optimum design of transport services, in terms of selecting service routes, scheduling frequencies, and allocating capacities, in order to meet customer requests at minimum cost. Classical SND models are typically formulated as deterministic problems, while in actual freight transport networks a high level of uncertainty exists. Travel time uncertainty, stochastic demand, limited fleet capacity, and other service disruptions all impact service reliability and operational efficiency. In multimodal networks, where services include more than one transportation mode such as road, rail, or waterways, the impact and complexity added by these uncertainties is further emphasized by the need to coordinate across operating conditions and service schedules. To address these challenges, recent research began incorporating stochastic elements into SND models, examining approaches such as scenario-based stochastic programming, robust optimization, and simulation-based approaches. However, there are still significant limitations. Many studies continue to focus on unimodal systems or treat multimodality in a simplistic manner. Furthermore, most existing approaches do not consider the option of modifying decisions dynamically in response to real-time disruptions. In particular, unification of tactical-level decisions (e.g., service selection and capacity planning) and operational-level decisions (e.g., truck allocation, vehicle routing, and dynamic re-planning) has received little attention.

To fill these gaps, this paper proposes a hybrid approach that integrates metaheuristics, simulation and machine learning techniques, which is able to address a complex version of the SND. The algorithm is designed to address the SND problem faced by an LSP operating in a multimodal freight transport network, where trucks operate flexibly while fixed train and barge schedules must be explicitly incorporated. Firstly, an algorithm is designed to solve the deterministic SND, which is later extended to address the problem with stochastic road travel times, finite truck fleet capacity, and adaptive re-planning decisions to face truck delays. The approach incorporates tactical and operational decision-making, enabling the generation of robust plans that can adjust dynamically to disruptions. A simulated annealing (SA) metaheuristic is designed to solve the SND obtaining high-quality solutions for large-scale problems. This is integrated with a simulation model, where operational complexities like dynamic truck allocation, vehicle routing, and re-planning policies due to travel time variability are included. To mitigate the computational burden of evaluating numerous candidate solutions, a regression-based surrogate model that estimates the costs is integrated into the optimization routine to allow for fast and adaptive decision-making. The numerical results highlight the effectiveness of the proposed methods in addressing both deterministic and stochastic problems. They also show that, despite relying on relatively simple techniques (SA and a regression model), their integration delivers strong performance without requiring more complex or sophisticated approaches. This demonstrates that practical, scalable solutions for real-world and large-scale freight transport planning can be achieved without excessive algorithmic complexity, making the approach both accessible and adaptable for logistics practitioners.

The remainder of this paper is organized as follows. Section 2 presents a review of the relevant literature. Section 3 provides a detailed description of the problem under investigation, describing the deterministic and stochastic versions of the SND. Section 4 outlines the proposed hybrid solution approach, describing the different components and its integration. Section 5 reports the experimental results and evaluates the performance of the proposed method. Finally, Section 6 concludes the paper and discusses potential directions for future research.

# 2 Literature review

This section reviews relevant literature for this study. Firstly, Section 2.1 examines research on service network design under uncertainty, with particular emphasis on multimodal systems, scheduled services, and stochastic travel times. Then, Section 2.2 provides a brief overview of how machine learning techniques have been used to enhance optimization methods in related freight transport planning problems. Finally, Section 2.3 highlights the main research gaps and positions the contributions of this paper within the existing literature.

## 2.1 Service network design under uncertainty

SND is a well-established problem that models the tactical decisions of an LSP in a freight transport network. Numerous versions of the problem exist, differing in assumptions, decision variables, and constraints. Comprehensive reviews of traditional models and solution approaches can be found in Elbert et al. [1], Crainic and Hewitt [2], and Crainic and Rei [3]. In its basic formulation, SND takes as input a graph representing the physical network and a set of transport orders or requests, sometimes also referred as commodities. The key decisions involve selecting which services to operate, determining their capacities or frequencies, and assigning cargo flows to these services. The objective is to meet transport demand while minimizing costs for the LSP. Variants of SND introduce extensions such as traffic flow balancing, crew planning, dynamic operations, and multimodal networks. Given its tactical nature, numerous studies have explored SND under uncertainty. This paper examines SND in a multimodal context with scheduled services and uncertain road travel times, so the rest of the literature review focus on those aspects, particularly restricted to studies that consider uncertainty

The SND under uncertainty has received less attention than its deterministic counterpart, though interest has grown in recent years. Delbart et al. [4] provide an overview of early studies in this area, particularly those considering intermodal or multimodal transport systems. They note that most studies focus on a single mode of transport, identifying only 14 papers between 2009 and 2018 that explicitly incorporate multiple transport modes. Since 2018, research has continued to concentrate primarily on unimodal systems. Our review identifies 17 relevant papers published since 2018 that address SND under uncertainty (see Table 1), of which only six consider multimodal systems. In general, the modes of transport correspond to a combination between road, rail or waterways [5–9], although alternatives such as hydrogen pipelines, combined with trucks, can also be considered [10]. Most recent studies consider demand as the main source of uncertainty, although some incorporate uncertainties related to transit times, available service capacities, or transport costs. Regarding objective functions, most models minimize expected costs, incorporating elements such as transit, storage, and transfer costs, as well as delay penalties and pollutant emissions quantifications. Ghavamifar et al. [11] use an alternative objective, as they maximize profits by also optimizing price levels. Ghanei et al. [12] adopt a multi-objective approach that seeks to maximize system fairness alongside cost minimization. Decision variables also vary substantially across studies. Most models include decisions on freight flow distribution and service selection or scheduling, with several extending to additional decisions from related problems such as facility location [13,14] and transport outsourcing [15,16]. Furthermore, some models incorporate fleet management aspects, including fleet allocation [15] and vehicle routing [7,15,17].

**Table 1: Recent literature on SND with uncertainties.**

| Reference | Uncertainty source | | | | Objective | Problem characteristics | | | | Method |
|---|---|---|---|---|---|---|---|---|---|---|
| | DL | TT | SC | TC | | MM | FM | TD | RR | |
| Ghavamifar et al. [11] | | | x | | profit | | | | | Exact |
| Hrušovský et al. [5] | | x | | | cost | x | | x | | Exact / simulation-based optimization |
| Layeb et al. [6] | x | x | | | cost | x | | x | | Simulation-based optimization |
| Richter and Stiller [19] | x | | | | cost | | | | | Heuristic |
| Zhao et al. [9] | x | x | | | cost | x | | x | | Metaheuristic |
| Hewitt et al. [15] | x | | | | cost | | x | x | | Matheuristic |
| Saedinia et al. [14] | | | | x | cost | | | | | Exact |
| Uddin and Huynh [8] | | | x | | cost | x | | x | | Exact |
| Wang and Qi [18] | x | | | | cost | | | x | | Exact |
| Baubaid et al. [20] | x | | | | cost | | | | | Exact / heuristic |
| Govindan and Gholizadeh [13] | x | | | | cost | | | | | Metaheuristic |
| Lanza et al. [21] | | x | | | cost | | | x | | Exact / metaheuristic |
| Müller et al. [7] | | x | | | cost | x | x | x | | Metaheuristic |
| Xiang et al. [16] | x | | | | cost | | | x | | Exact |
| Deng et al. [10] | x | | x | | cost | x | | | | Exact |
| Ghanei et al. [12] | | | x | | cost, fairness | | | | | Exact |
| Liu et al. [17] | x | | | | cost | | x | | | Matheuristic |
| This paper | | x | | | profit | x | x | x | x | Metaheuristic / simulation-based optimization |

**Abbreviations used in the table:**
**-Uncertainty sources: demand level (DL), travel times (TT), service capacity (SC), transport costs (TC).**
**-Problem characteristics: multimodal (MM), fleet management (FM), time dependent model (TD), real-time replanning (RR).**

Most studies formulate SND under uncertainty as a stochastic mixed-integer programming model, assuming some parameters follow a known probability distribution. These are typically two-stage problems: the first stage involves tactical decision-making, while the second stage addresses operational decisions once uncertainties realize. For instance, excess demand may require rerouting or outsourcing, while transit times longer than expected may lead to schedule adjustments. These models commonly use scenario-based approaches to deal with the uncertainties. Alternative modelling approaches include chance-constrained models [8,9], robust optimization [14,18] and simulation-based methods [5,6]. While many of these studies propose models that are dynamic in the sense that they account for time-dependent decisions, none of them explicitly incorporate real-time replanning into the problem formulation. Although two-stage models incorporate, to some degree, reaction decisions to face uncertainty, the scenario-based approaches mean that these reactions are precomputed rather than dynamically adapted.

Therefore, they develop plans that are optimal across a limited set of predefined scenarios, instead of modelling truly dynamic decision-making in a real-time context.

In terms of solution methods, most studies propose exact methods, either using commercial solvers to solve the formulated model [8,12,14,21] or developing decomposition and column-generation techniques [10,11,16,18]. Due to the problem complexity, heuristic approaches have also been explored, including ad-hoc heuristics [19,20], matheuristics [15,17], and metaheuristics such as ant colony optimization [9], cross-entropy algorithm [13], and iterated local search [7]. In this paper we propose a simulation-based optimization technique, and there are studies in the recent literature that apply related techniques to the SND under uncertainty. Hrušovský et al. [5] propose an iterative approach in which a deterministic model generates a transport plan for a multimodal network, which is then tested for reliability using a simulation model. If the plan is found to be unreliable, the deterministic model is re-solved with additional constraints to avoid the identified unreliable connections. The simulation model accounts for uncertain travel times, leading to potential delays and missed connections. However, in cases of missed connections, freight is rerouted via trucks without additional re-planning decisions. Furthermore, the model does not incorporate fleet management or vehicle routing. Layeb et al. [6] adopt a different simulation-based optimization approach, integrating an optimization module addressed by a commercial optimization software, with a simulation module for the evaluation of candidate solutions. They address the problem of scheduling a set of services in a multimodal network, under uncertain demand and travel times. However, similar to the previous case, the simulation model does not incorporate fleet management or vehicle routing decisions, nor does it include re-planning mechanisms to mitigate uncertainties.

## 2.2 Machine learning in freight transport planning

The intersection of machine learning (ML) and combinatorial optimization has received growing attention in recent years, with surveys highlighting the advancements and opportunities of integrating ML with different optimization techniques [22,23]. This growing interest is also evident in research on planning and optimization within freight transport problems. However, the application of ML to solve the SND problem remains relatively limited. Most existing research has focused on classical combinatorial problems such as the Vehicle Routing Problem (VRP) and the Traveling Salesman Problem (TSP) [24–26]. These problems, fundamental to logistics and supply chain management, aim to determine cost-efficient vehicle routes while satisfying constraints like vehicle capacity and delivery time windows. As these can be considered subproblems within the broader SND framework, insights gained from VRP and TSP studies offer valuable guidance for developing new approaches to SND.

Recent studies have highlighted the efficacy of ML-based methods in tackling optimization challenges, in some cases surpassing traditional optimization techniques. Many propose end-to-end approaches, leveraging neural networks and reinforcement learning to generate complete solutions from scratch. While these methods can produce high-quality solutions rapidly, their training is typically time-consuming, as they require large sets of training data. Moreover, their ability to generalize to more complex problem definitions, such as those involving more decisions beyond vehicle routing, is often limited. Hybrid methods represent an alternative approach that can partially overcome these challenges. ML models can be embedded within heuristic or metaheuristic algorithms to enhance their efficiency, for instance by learning patterns of high-quality solutions to guide the search process or by selecting suitable operators during optimization. This way, hybrid approaches can combine the domain knowledge of optimization techniques with the adaptability of ML, leading to improved performance in solving complex problems [26,27]. In addition to hybrid strategies, other data-driven methods have been increasingly adopted. Techniques like surrogate modeling [28] and Bayesian optimization [29] approximate complex or computationally expensive cost functions, enabling faster solution evaluations. Similarly, inverse optimization has been employed to infer model parameters from observed behaviors, facilitating the calibration of models based on real-world data [30]. Other approaches can be used to aid traditional optimization algorithms and improve their performance, for example, by reducing the problem size through the prediction of variables that can be fixed [31–33], or by using ML to select and configure the best algorithm setting to solve a specific instance of the problem [34]. Furthermore, beyond direct optimization, ML applications also extend to various supporting aspects of transportation planning. Techniques such as deep neural networks and reinforcement learning are used to predict routing behavior and traffic flow, providing critical inputs for optimization models. Additionally, these methods enable demand forecasting, travel time prediction, preventive maintenance scheduling, and anomaly detection in transport data, showcasing their versatility in enhancing multiple aspects of transportation systems [35,36]. Together, these advancements illustrate the growing potential of integrating learning-based and data-driven techniques with established optimization methods, offering more flexible, scalable, and adaptive tools for addressing the complexities of modern transportation systems.

### 2.3 Contributions to the literature

The insights from this literature review highlight that, despite growing interest in stochastic formulations of the SND problem, significant gaps remain, particularly in multimodal transport systems and the development of adaptive solution methods. Most approaches rely on scenario-based stochastic programming or robust optimization, which require predefined uncertainty representations and offer limited flexibility for real-time decision-making. While simulation-based methods have been explored, they have primarily been used to evaluate the feasibility or resulting costs of candidate solutions, rather than actively integrating optimization and decision-making within the simulation framework. Additionally, key operational elements such as fleet management, vehicle routing, and adaptive re-planning remain underexplored. In contrast, the approach proposed in this paper advances the field by integrating a simulation model that not only accounts for uncertainty but also incorporates dynamic decision-making and the cascading effects of uncertain travel times. Unlike existing methods, our simulation-based optimization approach allows for real-time adaptations of the decisions within the simulation itself, resulting in a more flexible and adaptive framework. Moreover, the proposed framework allows to consider more realistic assumptions about the dynamics within the system, providing a more realistic representation of the problem. Finally, adaptive learning techniques are incorporated within the solution framework to overcome the computational cost of simulation, an aspect that also has not been included in the research of the SND problem. Overall, the framework integrates well-established yet flexible and scalable modeling and optimization techniques, which in combination provide an efficient and high-performing method for tackling a complex representation of the stochastic SND.

## 3 Problem description

This paper addresses the SND, where an LSP makes tactical decisions regarding transport request acceptance, capacity booking on train and barge services, and the allocation and routing of trucks. We consider two versions of this problem: deterministic and stochastic. The deterministic version follows the definition in Zhang et al. [37] and serves as a benchmark for comparison with existing literature. The stochastic version extends the deterministic model by incorporating uncertainty in road travel times and a limited truck fleet, which increases the operational impact of truck delays, enabling a more robust decision making at the tactical level.

### 3.1 Deterministic problem

We analyze an LSP operating in a multimodal long-haul logistics network. In this tactical planning problem, the LSP must decide which transport requests to accept and fulfil. Additionally, a set of scheduled train and barge services with predefined capacities is available, requiring the LSP to determine how many capacity slots to book for each leg in each service. The LSP also operates a fleet of trucks, which, in the deterministic problem, is assumed to be unlimited, so truck connections are always available and uncapacitated. The LSP also designs a transport plan, specifying for each transport requests the transport(s) mode(s) to be used and the route to be followed.

**Table 2: Notation on problem formulation.**

| Notation | Description |
|---|---|
| **Sets:** | |
| $G(N,A)$ | Graph representing physical network. |
| $R$ | Transport requests indexed by $r$. |
| $S$ | Train and barge services indexed by $s$. |
| $L$ | Service legs indexed by $l$. |
| $L_s \subseteq L$ | Set of legs comprising service $s$. |
| $P$ | Feasible paths indexed by $p$. |
| $P_r \subseteq P$ | Feasible paths for request $r$. |
| $P_l \subseteq P$ | Feasible paths that use leg $l$. |
| **Parameters**: | |
| $b_r$ | Reward of accepting request $r$ [EUR]. |
| $d_r$ | Size of request $r$ [containers]. |
| $c_p^{transit}$ | Transit cost of path $p$ [EUR/container]. |
| $c_p^{transfer}$ | Transfer cost of path $p$ [EUR/container]. |
| $c_p^{store}$ | Storage cost of path $p$ [EUR/container]. |
| $c_p^{delay}$ | Delay cost of path $p$ [EUR/container]. |
| $e_l$ | Available capacity in service leg $l$ [containers]. |
| $g_l$ | Cost of booking in service leg $l$ [EUR/container]. |
| **Decision variables:** | |
| $\mathrm{x_r}$ | Binary variable, 1 if request $r$ is selected. |
| $y_l$ | Number of slots (containers) booked in service leg $l$. |
| $z_{r,p}$ | Number of containers from request $r$ using path $p$. |

Table 2 presents the notation used in the problem formulation of this problem. The problem is formulated assuming a physical network $G(N,A)$, where $N$ is the set of nodes, which can correspond to the origin-destination locations of the transport requests and/or intermodal terminals, and $A$ is the set of arcs connecting them. In this study we assume $G$ is a complete graph, since all locations are connected by road with a known distance. Additionally, some pair of nodes can be connected by train or barge services (multiple services can connect the same OD pair at different times). Over this network, a path is defined as a feasible route, potentially consisting of multiple legs via truck, train, and/or barge. The feasibility of a path is determined by the requests' time window and the departure and arrival times associated with each leg. For truck legs, departure times are flexible and travel times depend on the distance traveled. In contrast, train and barge legs follow the predefined schedules of the corresponding services, with fixed departure and arrival times. Feasibility therefore requires temporal consistency along the entire path: each leg must be reachable given the arrival time from the previous leg, allowing sufficient transfer time at intermediate terminals. In addition, the resulting path timing must comply with the request's time window. The lower bound of the time window is enforced as a hard constraint, while the upper bound is treated as a soft constraint, with penalties incurred for late deliveries. The identification of feasible transport paths is a complex problem in itself, which for simplicity is implicit in the presented formulation. These feasible paths can be pre-computed using a method as the one described in Section 4.2. It must be noted that the composition of the feasible path set used may affect the performance of the solution algorithms for this problem. Under these assumptions, the decision variables of the problem are: the selection of transport requests to be served ($\mathrm{x_r}$), the number of capacity slots booked on each service leg ($y_l$), and the assignment of transport paths to each selected request ($z_{r,p}$).

$$Max\,Z = \sum_{r \in R}\left(b_r x_r - \sum_{p \in P_r} z_{r,p}\left(c_p^{transit} + c_p^{transfer} + c_p^{store} + c_p^{delay}\right)\right) - \sum_{l \in L} g_l y_l \qquad (1)$$

$$\sum_{p \in P_r} z_{r,p} = d_r x_r \qquad ,\forall r \in R \qquad (2)$$

$$\sum_{p \in P_l}\sum_{r \in R} z_{r,p} \le y_l \qquad ,\forall l \in L \qquad (3)$$

$$y_l \le e_l \qquad ,\forall l \in L \qquad (4)$$

$$\mathrm{x_r} \in \{0,1\} \qquad ,\forall r \in R \qquad (5)$$

$$y_l \in \mathbb{N}_0 \qquad ,\forall l \in L \qquad (6)$$

$$z_{r,p} \in \mathbb{N}_0 \qquad ,\forall r \in R\,,\forall p \in P_r \qquad (7)$$

The objective function (1) aims to maximize the LSP's profit. Revenue is generated from the rewards associated with selected transport requests, while costs include transportation expenses (in-vehicle transit, transfers, storage, and delay penalties) and capacity booking costs for train and barge services. Constraints (2) ensure that a transport plan is defined for all accepted requests, allowing for the splitting of a request across multiple paths. Constraints (3) guarantee that the booked capacity is sufficient to transport the assigned containers for each service leg, while constraints (4) restrict the total booked capacity to the maximum available for each service leg. Constraints (5)-(7) specify the binary and integer nature of the decision variables.

## 3.2 Stochastic problem

The stochastic version extends the deterministic problem by incorporating uncertain road travel times and a limited truck fleet. These assumptions affect both the problem definition and the modelling approach. In particular, this problem emphasizes the relationship between tactical and operational level decisions and the feedback loop between them in order to adapt the decisions accordingly. A schematic view of the main characteristics of the stochastic problem is provided in Figure 1, highlighting the main decisions at each decision level. It must be noted that operational re-planning decisions are taken dynamically in a rolling-horizon manner as uncertainty unfolds. Therefore, the problem does not correspond to a classical two-stage stochastic programming model. The following subsections describe in more detail the specific assumptions and implications when incorporating uncertain travel times and the limited truck fleet in the problem.

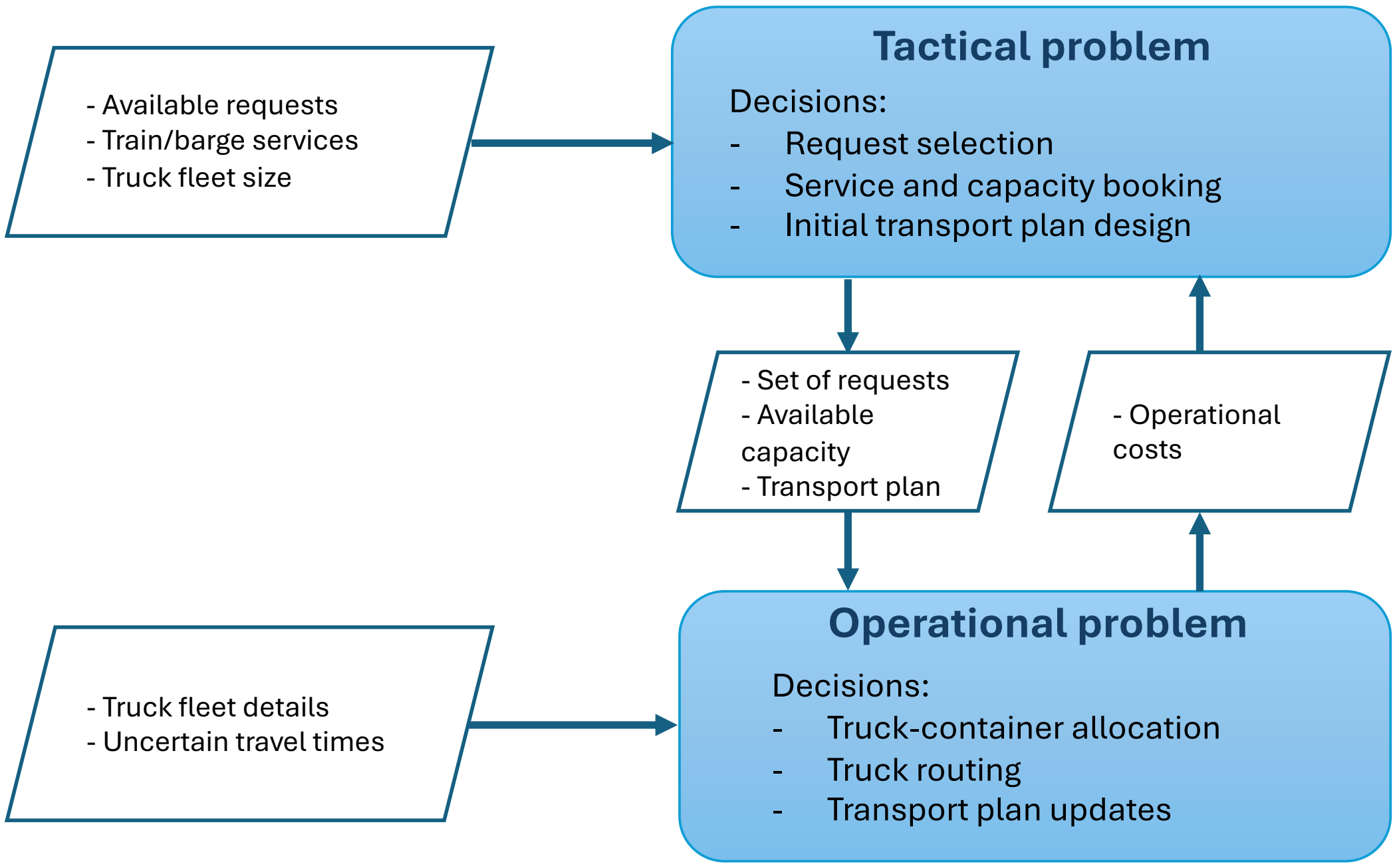

**Figure 1: Interaction between tactical and operational level decisions in the stochastic problem.**

### 3.2.1 Uncertain travel times

We assume that road travel times are uncertain due to two main factors: road disruptions and individual truck variability. Road disruptions introduce systematic delays affecting all trucks traveling on a given link for a certain period. These are events that normally happen with a low probability and may have a big impact on the affected truck trips. Additionally, travel times may vary between individual trucks, introducing minor variability compared to that caused by disruptions, but that affect all truck trips during the simulation. Then, the road travel time for a given truck at a given occurrence is determined by incorporating both sources of uncertainty into the baseline expected travel time, as shown in Equation (8), where $\tilde{t}_{ij}$ represents an instance of the actual travel time between locations $i$ and $j$, $\tilde{\eta}_{ij}$ accounts for ongoing disruptions that systematically affect all trucks on the road segment, $\tilde{\varepsilon}$ accounts for the individual variations in travel time for each truck, and $\bar{t}_{ij}$ denotes the baseline expected travel time in the absence of disruptions or individual variations.

$$\tilde{t}_{ij} = \left(1 + \tilde{\eta}_{ij}\right)(1 + \tilde{\varepsilon})\,\bar{t}_{ij}, \quad \tilde{\eta}_{ij} \in [0, \eta_{\max}], \quad \tilde{\varepsilon} \in [\varepsilon_{min}, \varepsilon_{max}] \tag{8}$$

Besides the particular trip between two nodes, these uncertainties can cause cascading delays throughout the transport network. For instance, a delayed truck may fail to meet a scheduled connection with a train or barge, resulting in a significantly larger overall delay for the corresponding request. Moreover, if the truck is assigned to transport multiple requests sequentially, its delay in one task may propagate to subsequent assignments. Therefore, delays of trucks may trigger the need to adapt the original travel plans. To handle these uncertainties, the problem is structured as a bi-level decision process. At the tactical level, the LSP selects transport requests, books capacity on train and barge services, and drafts an initial transport plan. However, as disruptions occur or delays propagates through the plan, this plan must be dynamically adapted during operations to accommodate these variations in travel times.

### 3.2.2 Limited truck fleet

In contrast to the deterministic version, the stochastic problem considers a finite truck fleet, requiring the explicit allocation of containers to specific trucks and the determination of their routing sequences, specifying the order in which various transport requests will be served. It is assumed that the LSP possesses a certain number of trucks, each assigned to a designated depot. Consequently, each truck's route must commence and conclude at its respective depot within the planning horizon. The truck fleet introduces additional planning requirements, thereby increasing the complexity of the problem. Furthermore, as previously noted, delays in the transportation of one request may affect subsequent requests scheduled on the same truck, increasing the propagation of delays. To address this, it is assumed that the truck routing plan is flexible, allowing for adjustments to routes and container allocations based on the system's current state. Accordingly, a provisional initial travel plan is generated at the tactical level, outlining the planned route for each request. Nevertheless, this plan is likely to be revised, meaning the actual routes followed, and consequently the associated cost components, are only fully known at the operational level.

### 3.2.3 Objective function

The characteristics of the stochastic problem imply that not all costs can be fully determined in advance, and additional unforeseen costs may arise during operations, which must be assumed by the LSP. In this context, the objective function for the stochastic problem maintains a structure similar to that of the deterministic case but considering the expected cost for the uncertain components. However, unlike classical stochastic optimization approaches, the uncertain costs in this problem do not follow a known probability distribution. Rather, they emerge from complex interactions between uncertain travel times and routing decisions within a multimodal transport network. Given the flexible routing approach assumed in this problem, relevant additional costs at the operational level may be related to extra kilometers driven, extended storage durations, or additional transshipments. Furthermore, delays may result in late deliveries to the LSP's customers, potentially incurring long term reputational damage in addition to potential fines or compensations. To account for this, a penalty factor is incorporated, proportional to the number of hours a container arrives late at its destination. Given these considerations, the LSP's profit maximization objective is formulated in Equation (9). The revenue from accepted requests ($\mathcal{R}$) and the cost of the booked capacity in scheduled services ($C_{cap}$) are defined at the tactical level. The expected value is considered for the cost components that are uncertain, namely: in-vehicle transit ($\tilde{C}_{transit}$), transfer ($\tilde{C}_{transfer}$), storage ($\tilde{C}_{store}$), and delay penalizations ($\tilde{C}_{delay}$).

$$Z = \mathcal{R} - \left( C_{cap} + E\left(\tilde{C}_{transit}\right) + E\left(\tilde{C}_{transfer}\right) + E\left(\tilde{C}_{store}\right) + E\left(\tilde{C}_{delay}\right) \right) \quad (9)$$

# 4 Solution approach

Due to the complexity of the problem, heuristic procedures are developed to ensure reasonable computational times. An SA algorithm is proposed, in which the SA operators determine the main tactical decisions, specifically the selection of transport requests and the scheduled services to be used. Once these key decisions are made, container and truck routing decisions are handled using methods that differ between the deterministic and stochastic versions of the problem. Notably, for the stochastic version, a simulation model is incorporated to assess the impact of uncertain travel times and limited truck availability. To mitigate the computational burden of the simulation process, a simple surrogate model is proposed and embedded within the SA to estimate the costs. Moreover, a learning-based adaptive method is designed to improve the performance of this surrogate model. The following sections provide a detailed explanation of each component of the solution approach.

## 4.1 Simulated annealing algorithm

SA is a well-known metaheuristic algorithm that has been successfully applied to various optimization problems. In this study, SA is chosen for its simplicity and promising preliminary results when applied to the problem under investigation. In essence, SA corresponds to a single-solution iterative algorithm, in which at each iteration a neighbor solution is generated, and it is evaluated against the current one, selecting with a probabilistic criterion which one passes to the next iteration. There are multiple variants in specific components of SA, and in this paper, we use an algorithm with an adaptive cooling schedule and re-heating mechanism. In our implementation, the solutions are represented by vectors which determine which transport requests are selected and how much capacity is booked in each leg of the scheduled services. Under this design, the evaluation step is not trivial and varies according to the problem version. For the deterministic problem, the evaluation procedure includes a heuristic that designs a transport plan, determining the paths for each selected request using the selected services (see section 4.2), from which the total profits can be computed. For the stochastic problem, the evaluation procedure is enhanced with a surrogate model that is trained using the insights of a simulation model, to account for the uncertainties and additional decisions included in the problem (see section 4.3).

## 4.2 Deterministic problem evaluation

For the deterministic problem, the evaluation step of the SA is addressed using a heuristic method, which defines the path for each selected request using the booked capacity in scheduled services, to then calculate the profits. These paths account for request time windows, service schedules, and travel times by truck, if applicable. Feasible paths can be multimodal, concatenating legs with different transport modes and services, although simplifying restrictions are imposed. In particular, trucks may be used for first- and last-mile legs when necessary, with up to two intermediate train or barge service legs (i.e., a maximum of four legs in total). To improve computational efficiency, all feasible paths for each transport request are precomputed prior to the execution of the SA algorithm, generating a pool of feasible paths per request. Then, during the evaluation of each candidate solution, this pool is filtered according to the services selected and the corresponding capacities available in that solution.

The precomputed pool of feasible solutions is generated as follows. Firstly, for each pair of nodes (independent of demand), all the possible connections using only train and barge services are identified. At this stage, at most one transfer is allowed, ensuring sufficient transfer time where applicable. Next, for each transport request, feasible paths are constructed by considering its specific time window. These may include direct truck connections or combinations of truck and scheduled service legs through intermediate transfer nodes. At this stage, for first-mile truck legs, travel times are explicitly considered to ensure timely arrival for connection with scheduled services, if applicable. It is worth noting that the direct truck connection is always a feasible path, for every transport request. Then, the cost of each feasible path is calculated based on the cost components defined in Section 3.1.

A specific candidate solution can be then evaluated using this pool of feasible paths. This involves constructing the overall transport plan according to the procedure described in Algorithm 1, where $q_{r,l} \in Q$ is an auxiliary variable indicating whether request $r$ is assigned to a path including service leg $l$, and $P' \subseteq P$ denotes the filtered subset of precomputed feasible paths, where any path using a service leg without selected capacity in the candidate solution is discarded. The plan consists of the route to be followed by each container in each transport request, respecting the maximum capacities of the selected services. The algorithm begins by assigning all containers of each selected request $r \in R'$ to the cheapest available transport alternative, considering only the selected train and barge services (set $S'$). Next, service capacities are checked for each leg $l$ in each selected service, and if the load assigned to that leg ($q_l$) in the current plan exceeds the available capacity in that leg ($l_{cap}$), a reassignment procedure is triggered. In this step, some containers are rerouted to the next cheapest feasible path with available capacity. The request to be reassigned ($r^* \in R_l$, where $R_l$ is the set of requests currently assigned to a path that includes leg $l$) is selected as the one that minimizes the increase in total cost, by searching the corresponding next cheapest alternative path ($p^* \in P$). Then, the number of containers to be reassigned is calculated, considering current the excess load in the leg, the amount from request $r^*$ currently assigned to that leg ( $q_{r,l}$), and the available capacity in the new path considering $p^*$ the current transport plan. This reassignment process continues iteratively until all train and barge services respect their maximum capacity constraints. It is important to note that in this deterministic problem, the truck fleet is assumed to be unlimited, ensuring that the problem is always feasible. If necessary, any remaining demand can always be transported entirely by truck from the origin to destination, albeit at a potentially higher cost.

**Algorithm 1: Heuristic for deterministic solution evaluation.**

```
//Initialization
Q ← 0
for each r ∈ R'
        {Q} ← assign_cheapest_path(S', P', r)
end for
for each s ∈ S'
        for each l ∈ L_s
                while q_l > l_cap
                        //Search best alternative
                        for each r ∈ R_l
                                {r*, p*} ← get_next_cheapest_path(S', P', Q, r)
                        end for
                        //Determine and re-assign correct amount
                        delta ← min (l_cap, q_{r*,l}, p*.available_capacity(Q) )
                        {Q} ← re_assign(r*, p*, delta)
                end while
        end for
end for
```

## 4.3 Stochastic problem evaluation

The evaluation step of the stochastic problem is more complex because not all costs can be determined in advance. Additionally, unlike the deterministic case, the stochastic problem accounts for a limited truck fleet, making truck allocation and routing crucial considerations. Moreover, it assumes the LSP makes adaptive re-planning decisions to deal with the uncertain travel times. To address this complexity, a discrete event simulation model is developed to simulate the operations of the LSP over a given timeframe based on the tactical decisions previously described. A simulation-based approach is selected due to its flexibility in handling uncertainties and its ability to capture the cascading effects of delays and disruptions. The model also allows for dynamic planning decisions by the LSP, enabling them to adapt to disruptions and minimize their impact. The simulation provides an accurate estimation of expected costs, which can then be used within the SA framework to evaluate candidate solutions. However, running a full simulation for each candidate solution can be computationally expensive, particularly for large problem instances, making it unpractical for optimization purposes. To address this issue, this study proposes a simple surrogate model based on a regression function that estimates delay costs without the need to run the simulation for every candidate solution. Additionally, an adaptive method is introduced to enhance solution quality while keeping computational complexity manageable. The following subsections describe the simulation model and its integration with the SA algorithm.

### 4.3.1 Simulation model

The simulation model represents the operational level decisions of an LSP within a multimodal logistics network with uncertain road travel times and scheduled train and barge services. The main inputs to the model include the graph $G(N, A)$ representing the physical network, the scheduled train and barge services with reserved capacity, the size and locations of the truck fleet, the set of selected requests, the initial travel plan for the selected requests, and various other parameters defining the problem instance. As mentioned before, besides uncertain travel times, the stochastic model differs from the deterministic one by considering a limited fleet of trucks. This fleet is operated by the same LSP, who must decide on the allocation of containers to specific truck vehicles and determine the sequencing and routing of these trucks. Overall, two main decision processes are executed within the simulation, which in turn are addressed in the model by means of three heuristic methods. The two main decision processes are:

- Operationalize the initial transport plan: The initial transport plan (generated by Algorithm 1), specifies for each selected transport request the sequence of legs, indicating the origin, destination, transport mode, and the corresponding service for trains and barges. Before the start of the simulation, for each truck leg a specific truck must be assigned that can complete the task within the corresponding time window, which is determined by the request time window and the scheduled times of train and barge services used in subsequent legs. The truck-container allocation heuristic is called sequentially for each container and each truck leg, until all truck legs in the plan are assigned. However, it is possible that no truck in the fleet can complete a given truck leg within the time window. In such cases, the travel plan for that specific request is deemed infeasible, triggering a search for an alternative route using the container routing heuristic.
- Adapt transport plan after significant delay: During the simulation, the uncertain travel times can lead to truck delays, which are considered significant when planned routes become unfeasible (due to missed connections). When a truck arrives late to any node, it is checked if the intended time windows can be achieved, both for its current cargo as for later requests in its planned route. If the time window cannot be respected for a given truck leg, first it is checked if the truck leg can be assigned to a new truck, calling the truck-container heuristic. If that is not possible, then the request travel plan is scrapped, and a new route is searched using the container routing heuristic.

Therefore, the aforementioned decision processed are addressed by repeatedly calling simple heuristics, namely:

- Container-truck allocation: Containers are allocated to a specific truck in the fleet considering a truck leg (origin, destination, time window). The allocation is determined by searching for the best insertion, looping over the planned routes for every truck, selecting the one that lead to the smallest total travelled distance increment.
- Container routing: As mentioned before, the container routing heuristic is needed when the original travel plan for a request is deemed unfeasible, due to a train/barge connection that cannot be reached on time. In such case, a new route is searched using the precomputed pool of feasible paths (see Section 4.2), but in this case the feasibility of the truck legs is considered, which determined with the heuristic for the container-truck allocation. It is important to note that this method assumes a high degree of flexibility in the options to modify the transport plan, for example in the terms to changing the assignment of containers within scheduled services.
- Truck routing: During these decision processes, also the routing of the trucks is defined and adapted dynamically. As mentioned before, the route of the trucks is constructed sequentially, looking for the best insertion of truck legs when required.

It is worth noting that, instead of the aforementioned heuristics, the simulation model could incorporate alternative methods for truck allocation and routing optimization, as well as incorporating more complex interactions between multiple players, such as dividing decisions between LSPs and carriers. The aim of this paper is to propose an adaptive simulation-based framework that leverage insights from the simulation model for the optimization process, that in future work could be extended to incorporate such more complex modelling in the operational model.

In summary, in its current form, the simulation follows the process illustrated in Figure 2. Initially, the initial transport plan is completed (operationalized) by determining the container-truck allocation and, if necessary, updating container routes. The discrete event simulation then is executed, processing events such as requests becoming available for pickup, vehicle departures and arrivals, and the start and end of loading and unloading operations. The only uncertainty considered in the current setting is truck travel time, which is calculated for each truck leg according to the description given in section 3.2.1. When stochastic travel times cause a truck to arrive late at a destination, the travel plan is adapted if necessary, replanning the routing of trucks and containers as needed, using the aforementioned heuristics procedures. At the end of the simulation, the relevant outcome of the operations can be calculated. For the purposes of this paper, the relevant outcome value corresponds to the operation costs and delay penalizations.

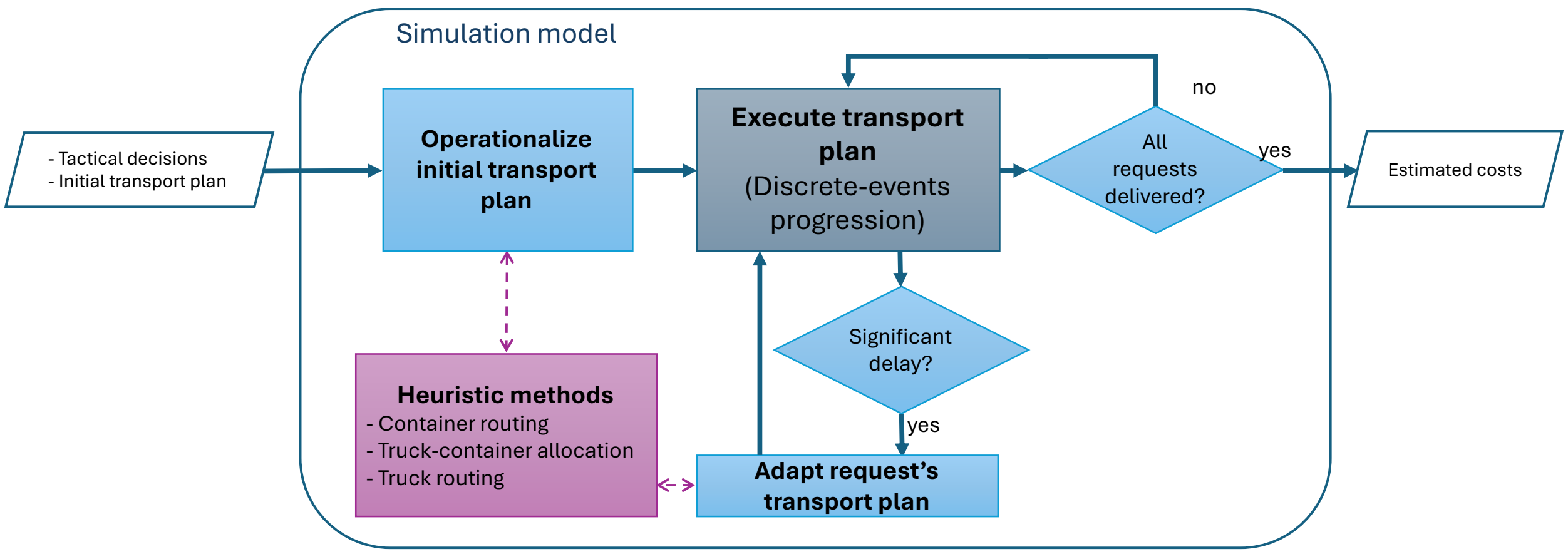

**Figure 2: Schematic view of the simulation process.**

### 4.3.2 Surrogate model to estimate costs

The simulation model described in the previous section allows for an accurate estimation of expected operational costs considering uncertain travel times and the replanning decisions made by the LSP at operational level. This estimation can be incorporated for more robust tactical level decisions of the LSP. However, the computation time of the simulation makes it impractical for optimization in large problem instances. To address this issue, a surrogate model is proposed to predict these costs for a given solution of the stochastic problem. Although multiple cost components vary at operational level, preliminary experiments showed that the most critical one corresponds to the delay penalizations. Therefore, the proposed model consists of a regression function that focuses on estimating these delay costs and penalizations. The proposed approach consists of firstly defining a scalar variable, denoted as $\gamma$, which aggregates relevant information about a candidate solution, and then fitting a function that takes $\gamma$ as input to estimate the total delay costs for that specific solution.

The variable $\gamma$ is defined so as to capture the key characteristics of the solution that influence delay costs, particularly related to the truck legs in the initial transport plan. For each of these legs, the number of containers to be transported, the expected truck travel time, and the time window for the transport leg are initially known. With this in mind, $\gamma$ is calculated by considering only the truck legs in the initial plan for each selected request. The number of assigned containers is multiplied by the expected travel time, including loading and unloading. These values are summed and divided by the available truck hours, which is the product of the available truck fleet size and the total time window for all the selected requests. The calculation follows equation (10), where *R* is the set of selected requests; $L_r$ is the set of truck legs in the initial plan for request *r*; $q_l$ is the number of container assigned to leg *l*; $t_l$ is the total expected time associated to leg *l*; *N* is the number of available trucks; and *W* is the length of the total time window considering the selected requests.

$$\gamma = \frac{\sum_{r \in R} \sum_{l \in L_r} q_l \cdot t_l}{N \cdot W} \tag{10}$$

A systematic correlation is found between γ and the delay costs obtained from the simulation model. To validate this relationship, a set of feasible solutions is generated for a given instance using the SA algorithm. To ensure diverse and representative solutions, selections are made from different stages of the SA process. Additionally, variations in problem parameters, such as the number of available transport requests and the relative size of the truck fleet, are introduced. For each feasible solution, the value of γ is calculated, and the simulation model is run to obtain expected delay costs. Figure 3 illustrates the relation between γ and the delay costs per container for the pool of solutions used in this study, showing a systematic correlation between both values, although it is a non-linear relation and at different scales. Therefore, a function is then fitted, using nonlinear least squares regression, to estimate the total delay costs using the value of γ. Various function forms are tested, with a cubic function yielding the best results. Then, equation (11) is used to estimate the total delay costs, and Figure 3 shows the comparison between the values of $f'(\gamma)$ and the total delay costs obtained with the simulation model. The results indicate that this method performs well even when problem parameters change, meaning the regression model can be pre-trained and applied across multiple problem instances. However, it must be noted that all the experiments in this paper were performed on instances with the same underlying physical network. Therefore, further testing is required to evaluate performance on different networks. If necessary, the function should be refitted with appropriate data to ensure reliability in varied scenarios.

$$\mathrm{E}(\mathrm{C}_{\mathrm{del}}) = f'(\gamma) = \mathrm{a}_0 + \mathrm{a}_1\gamma + \mathrm{a}_2\gamma^2 + \mathrm{a}_3\gamma^3 \tag{11}$$

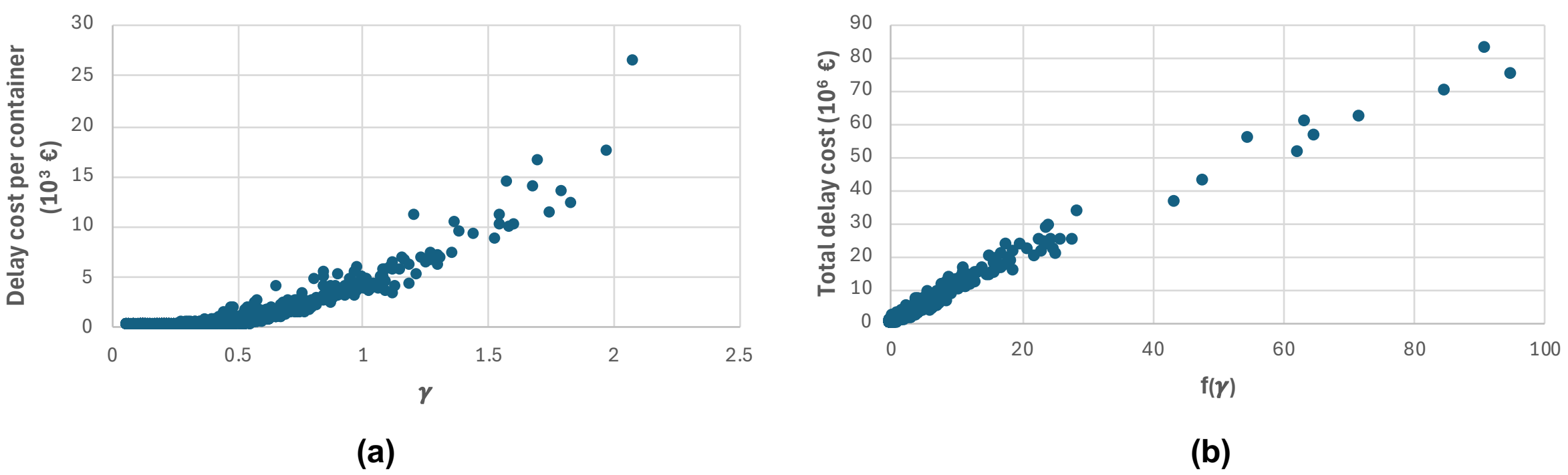

**Figure 3: Visualization of relation between γ and f′(γ) with simulation delay costs.**

### 4.3.3 Adaptive SA

The regression function proposed in the previous section can be incorporated into the SA framework to evaluate candidate solutions efficiently, allowing the generation of high-quality solutions without increasing computational time compared to deterministic methods. Since the function can be fitted once and then applied to multiple problem instances, it significantly reduces computational overhead. However, better results may be obtained if the function parameters are adjusted to the specific instance being solved. Instead of performing a complete refitting process for each instance, an adaptive method is proposed to dynamically adjust the parameters based on data collected during the SA execution.

The proposed adaptive method operates as follows. For most of the SA execution, candidate solutions are evaluated using the predefined regression function to estimate total delay costs. However, at every $n_1$ iterations, a set of $n_2$ candidate solutions, comprising the current solution and the neighbor solutions generated in previous iterations, is selected for evaluation using the simulation model, and the regression parameters are updated using the fitting procedure described earlier. Unlike the initial fitting process, which relies on a large dataset, this adaptation step uses only the $n_2$ selected solutions to update the regression parameters. To prevent excessive changes in the parameters, adjustments are capped at a maximum change of $\theta$%. Additionally, this adaptation is only applied after at least $n_3$ iterations have been completed from the start of the SA process, ensuring that modifications occur only when the solution quality has reached a relatively stable level. The parameters $n_1$, $n_2$, $n_3$, and $\theta$ must be carefully tested and tuned to balance the potential improvement in solution quality with the additional computational cost incurred by running simulation evaluations during SA execution (see Table 3 for our case).

## 4.4 Overview of the proposed SA variants

In summary, the method proposed in this paper correspond to the SA using the adaptive surrogate model, which in the following sections will be denoted $SA_A$. However, to test its performance, other variants of the SA are used in the experiments for benchmarking. Overall, five variants of the SA algorithm are defined, each differing in the evaluation procedure:

- $\boldsymbol{SA_H}$: Uses the heuristic procedure outlined in Algorithm 1 to evaluate each candidate solution. This serves as the baseline setting for solving the deterministic problem.
- $\boldsymbol{SA_B}$: Builds upon $SA_H$ by incorporating a 10% buffer to all road travel times when generating travel plans. Preliminary experiments indicated that this buffer helps mitigate delays caused by uncertain travel times. However, its exact value (around 10%) does not significantly affect the results. Therefore, this buffer remains fixed in the subsequent SA variants.

- $\boldsymbol{SA_F}$**:** Retains the heuristic evaluation approach from $SA_H$ but integrates the regression function as surrogate model (see Section 4.3.2) into the objective function for candidate solutions. As previously mentioned, this regression function was fitted once and applied across all experiments, based on a pool of feasible solutions encompassing instances of various sizes.
- $\boldsymbol{SA_A}$**:** Adaptive algorithm that extends $SA_F$ by dynamically adjusting the parameters of the regression function during the SA execution (see Section 4.3.3).
- $\boldsymbol{SA_S}$**:** Directly employs the simulation model to evaluate each candidate solution, yielding the most precise cost estimations. However, this approach comes at a significantly higher computational cost compared to the other variants.

In summary, the first two variants ($SA_H$ and $SA_B$) rely on fast heuristics for the evaluation step, but do not explicitly account for the additional complexities introduced in the stochastic problem. On the other hand, $SA_S$ uses the simulation model to evaluate candidate solutions, offering the true objective value for the stochastic problem. However, this approach is significantly more time-consuming. Therefore, $SA_F$ and $SA_A$ aim to bridge this gap, by using learning techniques to use the insights derived from the simulation model more efficiently. As such, $SA_F$ and $SA_A$ will be referred to as the learning-based algorithm in the subsequent sections.

# 5 Experimental results

We perform various experiments to evaluate the performance of the proposed methods. Firstly, we assess the results obtained by applying the SA algorithm to solve the deterministic problem, comparing these outcomes with results available in the literature. Secondly, since no previous studies have addressed the stochastic problem with the same characteristics, comparisons in this case are made among the various algorithmic settings proposed in this paper. The experiments are conducted on instances of varying sizes, though all are based on the same physical network and set of available scheduled services. This benchmark instance, originally introduced by Guo et al. [38], was also utilized by Zhang et al. [39] to evaluate multiple variants of the deterministic SND problem. The network consists of 10 terminals and 82 scheduled services, including both train and barge services. The size of the tested instances in this study is determined by the number of transport requests, which range from 5 to 1600. Further details regarding the benchmark instance can be found in Guo et al. [38]. The methods described in this paper were implemented in C++ and Python, and all experiments were performed on a MacBook Pro equipped with an Apple M3 Pro chip and 18 GB of RAM.

The experiments presented in this section employ the variants of the SA algorithm described in Section 4.1. While the general framework remains consistent across experiments, the evaluation procedure is adjusted based on the problem definition and the specific solution approach. Consequently, the primary algorithm parameters that define the SA are held constant throughout all experiments. Table 3 summarizes the values of these parameters, including those relevant to the adaptive version described in Section 4.3.3. These values were determined through preliminary experiments and by considering value ranges commonly adopted in the literature. Additionally, for each experiment, 30 replications of the SA procedure were performed, and average values are reported in the subsequent sections.

**Table 3: Selected values for SA parameters.**

| Parameter | Symbol | Value |
| --- | --- | --- |
| **Base SA:** | | |
| Maximum number iterations | – | 2000 |
| Number of iterations without improvement before reheating | – | 100 |
| Initial temperature | – | 1000 |
| Reheat temperature | – | 500 |
| Initial cooling rate | – | 0.99 |
| Interval for adaptive cooling rate | | (0.95, 0.999) |
| **Adaptive SA:** | | |
| Iterations between function update | $n_1$ | 100 |
| Solutions evaluated for function update | $n_2$ | 1 |
| Minimum iterations before function update | $n_3$ | 500 |
| Cap for parameters change | $\theta$ | 0.1 |

## 5.1 Benchmarking based on the deterministic problem

The first set of experiments is designed to evaluate the performance of the basic version of the algorithm ($SA_H$) on the deterministic version of the problem, for which we have detailed results in the literature. Specifically, we compare our approach to the state-of-the-art algorithm proposed by Zhang et al. [39], who introduced an Adaptive Large Neighborhood Search (ALNS) algorithm to solve large instances of the SND with flexible service options. For this comparison, we adopt the problem denoted as $L_1$ in Zhang et al. [39], which corresponds to a scenario involving flexible trucks alongside scheduled train and barge services. This setup aligns with the deterministic problem described in Section 3.1, with the additional conditions that all transport requests must be accepted and cannot be split across multiple routes. We modify our algorithm to fit this definition. Table 4 presents the comparison of results across instances ranging from 5 to 1600 transport requests. For each instance, the table reports the average cost of the best solutions obtained by the ALNS in Zhang et al. [39], alongside the best and average costs achieved by $SA_H$, as well as the relative difference between the two methods. The results indicate that $SA_H$ performs competitively, improving upon the solutions in most instances by a small margin, particularly in the larger ones. Additionally, a significant improvement in computational efficiency is observed. For example, in the largest instance involving 1600 transport requests, Zhang et al. [39] report requiring approximately 40 hours of computing time to obtain their best solutions. In contrast, $SA_H$ is able to find its reported solution in under one minute. This highlights not only the effectiveness of the SA framework in exploring the solution space, but also the efficiency of the heuristic procedure employed to solve the routing subproblem and evaluate solutions.

**Table 4: Comparison of results for the deterministic problem.**

| Instance (# of requests) | Total cost (EUR) | | | Diff (avg) |
|---|---|---|---|---|
| | Benchmark [39] | $SA_H$ (best) | $SA_H$ (avg) | |
| R-5 | 4269 | 4240 | 4262 | -0.2% |
| R-10 | 25 988 | 25 802 | 26 547 | 2.2% |
| R-20 | 42 776 | 42 375 | 43 435 | 1.5% |
| R-30 | 64 938 | 64 533 | 65 328 | 0.6% |
| R-50 | 131 556 | 125 905 | 126 742 | -3.7% |
| R-100 | 176 277 | 171 955 | 172 929 | -1.9% |
| R-200 | 480 447 | 467 916 | 470 642 | -2.0% |
| R-400 | 1 109 777 | 1 097 795 | 1 099 892 | -0.9% |
| R-700 | 1 064 366 | 1 046 949 | 1 048 025 | -1.5% |
| R-1000 | 1 021 288 | 999 049 | 1 000 174 | -2.1% |
| R-1300 | 1 053 422 | 1 030 712 | 1 031 501 | -2.1% |
| R-1600 | 1 031 410 | 1 017 702 | 1 018 497 | -1.3% |

## 5.2 Numerical experiments for the stochastic problem

Following the demonstration of the proposed solution approach's strong performance on the deterministic problem, the next step is to evaluate its effectiveness in addressing the more complex stochastic problem. However, as no prior studies have investigated this specific problem definition, there are no available results in the literature for direct comparison. Therefore, the experiments focus on testing and comparing the five algorithmic configurations proposed in Section 4.4, using the deterministic version ($SA_H$) as a baseline.

### 5.2.1 Experimental setting

The experiments are conducted on the same benchmark instances used in the deterministic case. Specifically, four instances, with sizes ranging from 50 to 400 transport requests, are considered. These instances are based on the dataset presented by Zhang et al. [39], from which the descriptions of the transport request sets are adopted. However, it was observed that the results on instance R-100 exhibited distinct trends compared to the others. To investigate the causes of this behavior, we made a new random selection of 100 transport requests from the overall set, which corresponds to the instance R-100' in the following experiments.

In addition, several problem elements are configured in accordance with the stochastic problem definition outlined in Section 3.2. In this sense, different scenarios are designed to capture variability in road travel times and variations in the size of the available truck fleet. Table 5 summarizes the scenarios considered for the stochastic problem. Two distinct scenarios are defined for each factor, which are then combined to yield four scenarios in total. Regarding travel time variability (see Section 3.2.1), the disruption component remains constant, modelled as a Poisson process with a mean interarrival time of 15 hours, durations uniformly distributed between 1–10 hours, and severities uniformly distributed with a maximum severity factor ($\eta_{\max}$) of 1. On the other hand, the individual truck variability is modeled using a beta distribution that bounds the variability relative to the baseline travel time. The two variability scenarios (low and high) are created by adjusting the upper bound of the individual variation component ($\varepsilon_{max}$), while $\varepsilon_{min}$ is fixed at -0.1. Concerning the truck fleet size, the fleet is arbitrarily distributed across the physical network, and the total number of trucks is determined relative to the instance size by multiplying the number of transport requests by a fleet factor $\kappa_{fleet}$, therefore defining two scenarios (small and large fleet). It must be noted that in these experiments a homogeneous fleet is assumed, although the proposed methods could address instances with trucks with varying costs or speeds. Among the four resulting scenarios, $\mathcal{V}^-\mathcal{F}^+$ represents a less constrained environment where the added complexity is expected to have minimal impact. Conversely, $\mathcal{V}^+\mathcal{F}^-$ reflects a more constrained setting, in which the influence of variability and limited fleet size is anticipated to be more significant.

**Table 5: Considered scenarios for the stochastic problem.**

| Scenario | Description | Parameter value | |
|---|---|---|---|
| | | $\varepsilon_{max}$ | $\kappa_{fleet}$ |
| $\mathcal{V}^-\mathcal{F}^+$ | low variability-large fleet | 0.1 | 0.5 |
| $\mathcal{V}^+\mathcal{F}^+$ | high variability-large fleet | 0.25 | 0.5 |
| $\mathcal{V}^-\mathcal{F}^-$ | low variability-small fleet | 0.1 | 0.25 |
| $\mathcal{V}^+\mathcal{F}^-$ | high variability-small fleet | 0.25 | 0.25 |

### 5.2.2 Comparison of SA variants

In line with the definitions presented in Section 4.4, a first set of experiments is conducted to evaluate the performance of five variants of the SA algorithm, namely: $\boldsymbol{SA_H}$, $\boldsymbol{SA_B}$, S$\boldsymbol{A_F}$, $\boldsymbol{SA_A}$, and $\boldsymbol{SA_S}$. For comparison purposes, the final best solution from each replication of the five SA settings is re-evaluated using the simulation model, and these values are reported. Moreover, each time the simulation model is employed, it is executed five times per solution, and the average cost results are used for the solutions objective function evaluation. Table 6 illustrates the comparative performance of the five SA settings across different instances and scenarios. For each combination of instance and scenario, the table presents the average profit (in euros) for the best solution found, as well as the average CPU time (in seconds), for each SA setting. From these results, multiple aspects can be discussed.

**Table 6: Comparison of results on multiple scenarios.**

| Algorithm | Profit (€) | | | | CPU time (s) | | | |
|---|---|---|---|---|---|---|---|---|
| | $\mathcal{V}^-\mathcal{F}^+$ | $\mathcal{V}^+\mathcal{F}^+$ | $\mathcal{V}^-\mathcal{F}^-$ | $\mathcal{V}^+\mathcal{F}^-$ | $\mathcal{V}^-\mathcal{F}^+$ | $\mathcal{V}^+\mathcal{F}^+$ | $\mathcal{V}^-\mathcal{F}^-$ | $\mathcal{V}^+\mathcal{F}^-$ |
| **R-50:** | | | | | | | | |
| $SA_H$ | 353 511 | 334 936 | 14 845 | -326 163 | 0.1 | 0.1 | 0.1 | 0.1 |
| $SA_B$ | 361 172 | 356 351 | 188 070 | 162 275 | 0.1 | 0.1 | 0.1 | 0.1 |
| $SA_F$ | 359 092 | 354 748 | 315 693 | 317 492 | 0.1 | 0.1 | 0.1 | 0.1 |
| $SA_A$ | 357 686 | 353 655 | 328 092 | 327 639 | 82.2 | 80.7 | 77.9 | 79.7 |
| $SA_S$ | 369 727 | 365 673 | 350 116 | 344 515 | 63.7 | 65.9 | 49.8 | 50.1 |
| **R-100:** | | | | | | | | |
| $SA_H$ | 827 862 | 824 779 | 812 202 | 792 695 | 0.2 | 0.2 | 0.2 | 0.2 |
| $SA_B$ | 825 871 | 824 274 | 813 271 | 814 437 | 0.2 | 0.2 | 0.2 | 0.2 |
| $SA_F$ | 824 568 | 823 694 | 803 723 | 807 084 | 0.2 | 0.2 | 0.2 | 0.1 |
| $SA_A$ | 812 878 | 809 316 | 812 021 | 814 236 | 79.9 | 78.5 | 79.4 | 78.3 |
| $SA_S$ | 825 626 | 824 596 | 821 581 | 821 112 | 102.4 | 106.6 | 94.7 | 95.5 |
| **R-100':** | | | | | | | | |
| $SA_H$ | 756 723 | 751 010 | 564 338 | 58 623 | 0.2 | 0.2 | 0.2 | 0.2 |
| $SA_B$ | 760 734 | 758 207 | 629 000 | 486 989 | 0.2 | 0.2 | 0.2 | 0.2 |
| $SA_F$ | 755 825 | 756 267 | 706 860 | 705 649 | 0.2 | 0.2 | 0.2 | 0.1 |
| $SA_A$ | 750 430 | 750 142 | 730 507 | 724 399 | 78.5 | 78.3 | 77.9 | 78.0 |
| $SA_S$ | 761 904 | 759 093 | 746 636 | 739 548 | 158.4 | 167.8 | 145.4 | 136.9 |
| **R-200:** | | | | | | | | |
| $SA_H$ | 1 527 763 | 1 505 518 | -45 686 | -1 247 808 | 0.5 | 0.5 | 0.5 | 0.5 |
| $SA_B$ | 1 522 342 | 1 521 475 | 632 515 | 530 289 | 0.5 | 0.5 | 0.5 | 0.5 |
| $SA_F$ | 1 502 628 | 1 502 931 | 1 348 591 | 1 341 008 | 0.5 | 0.4 | 0.3 | 0.3 |
| $SA_A$ | 1 498 240 | 1 502 468 | 1 369 742 | 1 382 958 | 78.3 | 81.6 | 78.1 | 79.3 |
| $SA_S$ | 1 522 030 | 1 521 432 | 1 471 319 | 1 459 091 | 793.6 | 806.3 | 539.9 | 572.8 |
| **R-400:** | | | | | | | | |
| $SA_H$ | 2 940 480 | 2 913 828 | -6 386 460 | -9 086 212 | 1.9 | 1.9 | 1.8 | 1.8 |
| $SA_B$ | 2 914 474 | 2 913 795 | -2 299 670 | -2 546 095 | 1.6 | 1.6 | 1.8 | 1.6 |
| $SA_F$ | 2 824 687 | 2 824 134 | 2 500 846 | 2 501 052 | 1.3 | 1.2 | 0.8 | 0.8 |
| $SA_A$ | 2 839 357 | 2 761 760 | 2 565 300 | 2 561 040 | 78.7 | 79.6 | 79.2 | 78.3 |
| $SA_S$ | 2 909 900 | 2 904 149 | 2 739 470 | 2 736 945 | 3057.6 | 3194.8 | 1857.1 | 1885.3 |

Firstly, it is possible to draw conclusions regarding the impact of the added complexity introduced in the stochastic model (uncertain travel times and truck routing decisions). In scenarios with a large truck fleet ($\mathcal{V}^-\mathcal{F}^+$ and $\mathcal{V}^+\mathcal{F}^+$), the impact appears to be minimal as expected since we do not have scarcity of resources. This is evidenced by the lack of significant improvement when comparing the performance of $SA_S$ to $SA_H$. Furthermore, in these scenarios, $SA_F$ and $SA_A$ tend to perform slightly worse than $SA_H$ or $SA_B$, suggesting that explicitly modelling this complexity does not consistently yield benefits. Conversely, in scenarios with a small truck fleet ($\mathcal{V}^-\mathcal{F}^-$ and $\mathcal{V}^+\mathcal{F}^-$), the difference in performance becomes pronounced, particularly as instance size increases. In such cases, the performance of $SA_H$ and $SA_B$ deteriorates substantially compared to $SA_S$. Indeed, for instance R-400, these simpler methods even yield solutions with negative profits. Regarding the different levels of travel time variability, no consistent significant effect is observed. Comparisons between $\mathcal{V}^-\mathcal{F}^+$ and $\mathcal{V}^+\mathcal{F}^+$, as well as $\mathcal{V}^-\mathcal{F}^-$ and $\mathcal{V}^+\mathcal{F}^-$, reveal minimal differences in outcomes. The only noticeable effect is the worse performance of $SA_H$ in $\mathcal{V}^-\mathcal{F}^-$ compared to $\mathcal{V}^+\mathcal{F}^-$. However, this difference is mitigated when $SA_B$ is applied, indicating that the buffer alone greatly mitigates the negative impact of additional variability. It is worth noting that preliminary experiments did show substantial effects under unrealistically high variability or severe disruptions, but these scenarios are outside the scope of interest for this study. The case of instance R-100 further highlights the importance of the relative fleet size. The results reveal no significant advantage of $SA_S$ over $SA_H$, even in small fleet scenarios (unlike all other instances). Upon closer examination of this instance's characteristics, it was found that, relative to other instances, the transport requests exhibited in average smaller sizes (in terms of the number of containers) and shorter distances between origin and destination. These features create a less constrained environment, reducing the demand for truck-hours relative to fleet capacity, with a similar effect to an increased truck fleet. Indeed, when making a new random selection of 100 transport requests (R-100'), whose characteristics are now in line with the instances of other sizes, the relative performance of the different algorithmic settings follow the same trends as for the other instances.

Figure 4 provides another perspective to compare the quality of the solutions, displaying the variability of the profits for the scenario $\mathcal{V}^-\mathcal{F}^-$, instead of just the average between replications. The plots show that $SA_H$ and $SA_B$, besides generating worst solutions in average, have much less consistent results. On the other hand, $SA_A$ presents a slightly higher variability between replications compared to $SA_F$ and $SA_S$, highlighting a potential aspect to improve in the adaptive method. Nevertheless, on these last three SA settings ($SA_F$, $SA_A$, and $SA_S$) the variability is greatly reduced, showing the consistency of the proposed algorithms. In conclusion, the advantage of applying the more complex methods proposed in this paper depends on the instance characteristics. Nonetheless, the significant performance gains observed in small fleet scenarios justify the additional effort.

The results presented in Table 6 also show the strong performance of the learning-based algorithms ($SA_F$ and $SA_A$), particularly in small fleet scenarios ($\mathcal{V}^-\mathcal{F}^-$ and $\mathcal{V}^+\mathcal{F}^-$). While they tend to perform slightly worse than $SA_S$ in terms of solution quality, the observed gap, up to 10% for $SA_F$ and 7% for $SA_A$, is outweighed by the computational efficiency. Notably, as instance sizes grow, the computational time for $SA_S$ increases rapidly, whereas it remains relatively stable and significantly lower for the learning-based algorithms. Of course, when comparing the CPU time of the learning-based algorithms it is relevant to also take into account the training time. As discussed in Section 4.3.2, the regression model was trained once for all the instances addressed in this paper. This was done with a pool of 19,500 solutions and the total training took 1113 s, of which most of the time corresponds to evaluate each solution with the simulation. Therefore, the training of the model took less computing time than one execution of $SA_S$ for the largest instance (R-400), and the difference would be larger if larger instances are addressed. Therefore, even if re-training was necessary to address instances based on a different physical network, for large instances the proposed learning-based method is expected to be faster, particularly if multiple experiments will be executed.

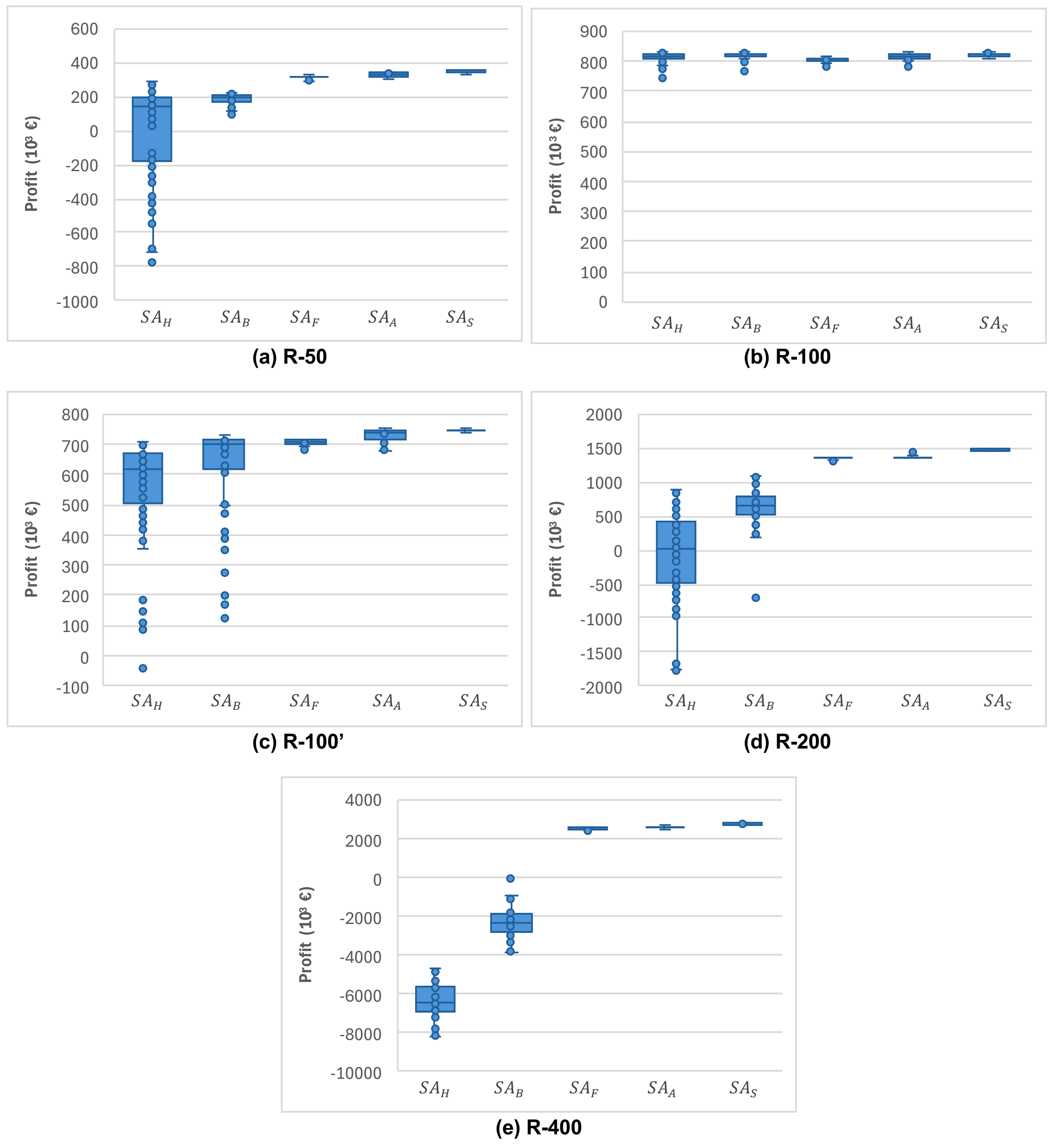

**Figure 4: Comparison of objective function value on scenario $\mathcal{V}^-\mathcal{F}^-$.**
***Note: The y-axis scale varies across plots to better reflect the range of values for each instance.***

**Table 7: Comparison of solution characteristics on scenario $\mathcal{V}^{-}\mathcal{F}^{-}$.**

| Algorithm | Profit components (x10³€) | | | | | | Selected requests | Mode share | | | Truck hours | Booked capacity | Used capacity |
|---|---|---|---|---|---|---|---|---|---|---|---|---|---|
| | Total | $\mathcal{R}$ | $C_{trs}$ | $C_{trf}$ | $C_{str}$ | $C_{del}$ | | $\sigma_R$ | $\sigma_M$ | $\sigma_S$ | | | |
| **R-50:** | | | | | | | | | | | | | |
| $SA_H$ | 14.8 | 478.7 | 64.2 | 37.8 | 6.1 | 355.7 | 0.96 | 0.07 | 0.38 | 0.55 | 1234 | 220 | 0.96 |
| $SA_B$ | 188.1 | 476.0 | 59.9 | 39.0 | 5.1 | 183.9 | 0.95 | 0.03 | 0.42 | 0.55 | 900 | 221 | 1.00 |
| $SA_F$ | 315.7 | 383.7 | 33.6 | 30.3 | 2.0 | 2.1 | 0.77 | 0.04 | 0.31 | 0.66 | 339 | 159 | 1.00 |
| $SA_A$ | 328.1 | 408.3 | 39.4 | 32.5 | 2.3 | 6.0 | 0.82 | 0.03 | 0.33 | 0.64 | 444 | 171 | 1.00 |
| $SA_S$ | 350.1 | 437.7 | 47.5 | 34.1 | 2.8 | 3.2 | 0.88 | 0.07 | 0.34 | 0.59 | 665 | 185 | 1.00 |
| **R-100:** | | | | | | | | | | | | | |
| $SA_H$ | 812.2 | 999.7 | 95.6 | 74.0 | 3.5 | 14.3 | 1.00 | 0.04 | 0.23 | 0.73 | 1377 | 406 | 0.98 |
| $SA_B$ | 813.3 | 998.3 | 97.4 | 74.4 | 3.2 | 10.0 | 1.00 | 0.03 | 0.23 | 0.74 | 1298 | 405 | 0.99 |
| $SA_F$ | 803.7 | 956.3 | 79.7 | 70.2 | 2.5 | 0.2 | 0.96 | 0.02 | 0.20 | 0.78 | 673 | 374 | 1.00 |
| $SA_A$ | 812.0 | 971.7 | 84.3 | 72.1 | 2.7 | 0.6 | 0.97 | 0.02 | 0.21 | 0.77 | 769 | 388 | 0.99 |
| $SA_S$ | 821.6 | 992.7 | 93.1 | 73.9 | 3.2 | 0.8 | 0.99 | 0.03 | 0.23 | 0.74 | 1093 | 401 | 0.99 |
| **R-100':** | | | | | | | | | | | | | |
| $SA_H$ | 564.3 | 980.0 | 145.3 | 70.5 | 5.6 | 194.2 | 0.98 | 0.10 | 0.29 | 0.61 | 2469 | 490 | 0.93 |
| $SA_B$ | 629.0 | 979.7 | 148.0 | 72.1 | 5.1 | 125.5 | 0.98 | 0.08 | 0.31 | 0.61 | 2219 | 488 | 0.95 |
| $SA_F$ | 706.9 | 876.7 | 99.6 | 65.5 | 3.5 | 1.2 | 0.88 | 0.02 | 0.25 | 0.73 | 770 | 429 | 1.00 |
| $SA_A$ | 730.5 | 920.0 | 113.2 | 68.4 | 3.5 | 4.5 | 0.92 | 0.05 | 0.26 | 0.69 | 1185 | 448 | 0.98 |
| $SA_S$ | 746.6 | 966.0 | 139.8 | 71.9 | 4.2 | 3.5 | 0.97 | 0.07 | 0.30 | 0.63 | 1821 | 484 | 0.96 |
| **R-200:** | | | | | | | | | | | | | |
| $SA_H$ | -45.7 | 1943.3 | 252.4 | 150.8 | 27.2 | 1558.6 | 0.97 | 0.08 | 0.39 | 0.53 | 5891 | 829 | 0.93 |
| $SA_B$ | 632.5 | 1934.3 | 248.7 | 154.1 | 23.2 | 875.8 | 0.97 | 0.06 | 0.41 | 0.54 | 4833 | 827 | 0.96 |
| $SA_F$ | 1348.6 | 1636.7 | 150.1 | 128.1 | 9.6 | 0.3 | 0.82 | 0.02 | 0.31 | 0.68 | 1450 | 686 | 0.99 |
| $SA_A$ | 1369.7 | 1680.7 | 159.0 | 132.5 | 11.2 | 8.3 | 0.84 | 0.02 | 0.33 | 0.65 | 1757 | 702 | 0.99 |
| $SA_S$ | 1471.3 | 1859.3 | 213.0 | 150.4 | 17.7 | 6.9 | 0.93 | 0.03 | 0.41 | 0.56 | 3403 | 790 | 0.99 |
| **R-400:** | | | | | | | | | | | | | |
| $SA_H$ | -6386.5 | 3842.7 | 612.4 | 277.3 | 92.3 | 9247.1 | 0.96 | 0.12 | 0.38 | 0.50 | 17110 | 1497 | 0.93 |
| $SA_B$ | -2299.7 | 3765.3 | 557.6 | 278.0 | 63.1 | 5166.3 | 0.94 | 0.10 | 0.38 | 0.52 | 12971 | 1474 | 0.96 |
| $SA_F$ | 2500.8 | 3026.0 | 282.8 | 225.9 | 16.2 | 0.3 | 0.76 | 0.05 | 0.28 | 0.67 | 3267 | 1237 | 1.00 |
| $SA_A$ | 2565.3 | 3142.7 | 308.0 | 236.6 | 18.0 | 14.7 | 0.79 | 0.05 | 0.31 | 0.64 | 4193 | 1277 | 1.00 |
| $SA_S$ | 2739.5 | 3457.7 | 421.0 | 263.7 | 26.3 | 7.2 | 0.86 | 0.07 | 0.38 | 0.55 | 7849 | 1375 | 0.99 |

**Notation: Total profits; revenue ($\mathcal{R}$); transit cost ($C_{trs}$) (including road transport train/barge services); transfer cost ($C_{trf}$); storage cost ($C_{str}$); delay penalties ($C_{del}$); proportion of available requests that are selected; proportion of requests transported solely by road ($\sigma_R$); proportion of requests transported only by train/barge scheduled services ($\sigma_S$); proportion of requests transported via a multimodal path ($\sigma_M$); total truck hours driven (both empty and loaded); total capacity booked in train/barge services (measured in container*km); ratio of booked capacity that is effectively used.**

Overall, these findings highlight the capability of the learning-based approaches to address more complex, constrained problem settings. However, it is important to analyze the performance gap between the proposed learning-based algorithms ($SA_F$ and $SA_A$) relative to $SA_S$ in greater detail, and this can be done by looking in more detail some characteristics of the solutions found by each SA setting. Table 7 displays some characteristics of the solutions found by each algorithm for scenario $\mathcal{V}^{-}\mathcal{F}^{-}$. By examining the cost components, it is evident that the poorer performance of $SA_H$ and $SA_B$ is largely attributable to high delay penalty costs. This results from their tendency to select more transport requests, which cannot subsequently be fulfilled on time due to the limited truck fleet, as the basic evaluation heuristic used in these methods does not account for fleet size constraints. In contrast, $SA_F$, $SA_A$, and $SA_S$ select fewer requests, successfully minimising delay penalties. However, $SA_F$ and $SA_A$ appear overly conservative, selecting fewer requests to such an extent that the income reduction is not compensated by the reduction in costs. This explains the performance gap relative to $SA_S$, and also why they may underperform $SA_H$ in large fleet scenarios ($\mathcal{V}^-\mathcal{F}^+$ and $\mathcal{V}^+\mathcal{F}^+$). Addressing this conservativeness signals a potential approach for improving the learning-based algorithms.

Table 7 further highlights other notable characteristics of the solutions produced by the different methods, such as the modal split. Due to the underlying physical network and scheduled services, the share of transport requests fulfilled solely by road ($\sigma_R$) is generally low, typically under 10%. However, this ratio is even lower In the solutions generated by $SA_F$ and $SA_A$, and to a lesser extent $SA_S$, and it is interesting to note that reduction is not compensated by multimodal transport ($\sigma_M$), but instead by paths using only train and barge scheduled services ($\sigma_S$). This shows how the proposed learning-based algorithms prioritize the selection of requests that can be served by these services. A similar trend is observed in the total truck hours, which are significantly lower in solutions generated by $SA_F$ and $SA_A$. Additionally, while these methods book slightly less capacity on scheduled train and barge services, the relative decrease is moderate. Importantly, $SA_F$, $SA_A$, and $SA_S$ yield solutions where the booked capacity is effectively utilised. In contrast, particularly in larger instances, $SA_H$ and $SA_B$ produce solutions where up to 7% of booked capacity remains unused due to missed connections and insufficient truck availability to implement recovery plans. These findings demonstrate that the proposed algorithms not only reduce overall costs and delay penalties but can also be used to obtain solutions with a more efficient use of the available capacities.

### 5.2.3 Managerial insights

The aforementioned results provide several practical insights for LSPs. Firstly, they highlight how relatively simple methods can be combined in an adaptive decision-making framework to efficiently address complex problems, such as the stochastic SND. Moreover, it is crucial to identify when it is necessary to explicitly model operational complexity and select the solution techniques accordingly. For example, in our problem definition, the use of the simulation models and the learning techniques seems only relevant when the truck availability is highly constrained, i.e. under scarce capacity. Similarly, it is likely that in extended problem definitions, that for example include more sources of uncertainty, more sophisticated metaheuristics or ML techniques will be required. Overall, selecting the right techniques aiming for fast solution methods is relevant even at tactical or strategic level, as they enable managers to run what-if scenarios. For example, our model assumes a high degree of flexibility in terms of the re-planning routing decisions to face disruptions, while in practice this is not always possible. Such experiments can be used to test different scenarios and quantify the potential gains and costs of adding extra flexibility and guide the efforts accordingly. Moreover, it is interesting to note that the algorithm solutions do not select all the available requests. Instead, the solutions reach a balance in the reward from accepting more requests and the associated risk of facing delays. The model can be used to explore in more detail this trade-off, and similarly to explore which routes to prioritize, both in terms of the selected demand and the services in which to book capacity. Finally, it is important to note that although in this work the simulation insights and the surrogate model were used to address the SND problem, this approach could be extended to other tactical problems with similar operational settings, providing managers with versatile tools to evaluate complex trade-offs and design more resilient transport plans.

# 6 Conclusions

This paper presents an adaptive simulation-based optimization framework to tackle the SND under uncertain travel times, addressing several key gaps identified in the literature. Specifically, the proposed approach integrates both tactical and operational decision-making processes, allowing for a realistic representation of cascading effects of delays and dynamic re-planning decisions. In particular, the model aims to maximize the profits for an LSP within a multimodal logistics network. An SA algorithm serves as the primary metaheuristic to efficiently explore the solution space at the tactical level, while the operational problem is modelled through a detailed discrete-event simulation. This simulation accounts for variability in road travel times due to disruptions and individual truck-level uncertainties, providing a more realistic assessment of operational outcomes. However, evaluating each candidate solution through full simulation incurs significant computational costs, especially for large-scale instances. To mitigate this, we introduce a regression function as surrogate model that predicts delay penalties based on key attributes of candidate solutions. Furthermore, an adaptive mechanism is incorporated within the SA framework, enabling the periodic refinement of the regression function based on simulation results collected during the optimization process.

The experimental results show the effectiveness of the proposed approach. For the deterministic version of the problem, the SA algorithm yields competitive results, outperforming existing benchmarks both in solution quality and computational efficiency. For the stochastic problem, the learning-based methods provide a good balance between solution quality and computational time. The advantages are particularly clear in instances with more constrained truck fleets, i.e., when there is scarcity in capacity. Nevertheless, results reveal less benefit in less constrained scenarios, i.e., with a surplus of resources, where the deterministic version of the algorithm already performs good. In fact, the conservative bias of the regression function in the selection of requests may lead to slightly worse results. This suggests a way to improve the methods by considering more details of the instance characteristics. Overall, the results demonstrate that although the used techniques are relatively simple in their respective domains (optimization and ML), the adaptive framework used to integrate them provides a robust method to tackle efficiently this complex version of the stochastic SND. Moreover, the computational overhead associated with training the learning component is modest relative to the solution time of large instances, and the model can be reused across instances defined on the same physical network. Therefore, even if re-training is required for instances based on different networks, the approach remains practically scalable.

There are several future research paths to improve the proposed methods and its applicability. Firstly, the problem definition can be extended to incorporate additional sources of uncertainty, such as fluctuating demand or service capacities, providing a more realistic and comprehensive modeling framework and analysis. Secondly, the simulation model could be expanded to include the interactions and negotiation dynamics among multiple logistics stakeholders, such as shippers, carriers, and LSPs in a multi-agent setting. Exploring these aspects would contribute to a richer understanding of real-world logistics systems and further improve the applicability of the proposed framework. Thirdly, more advanced learning-based techniques, such as neural networks or inverse optimization methods, could be investigated and incorporated as surrogate models, to further enhance the accuracy of cost estimation without incurring excessive computation. The results presented in this paper demonstrate that the proposed algorithm, including the simple regression model, performs well and provides high-quality solutions. Therefore, we consider such complex learning techniques unnecessary for the current problem formulation. Nonetheless, they may become valuable for extended problem variants, such as those incorporating additional uncertainties or multi-agent interactions.

In summary, this paper demonstrates that the integration of simulation models with metaheuristic optimization and learning techniques provides a powerful and adaptive solution framework for addressing complex SND problems. The results highlight the significant potential of such approaches in developing efficient decision-support tools for real-world logistics operations characterized by uncertainty and operational complexity.

## Acknowledgment

Funded by the European Union (ERC, ADAPT-OR, No 101117675). Views and opinions expressed are however those of the author(s) only and do not necessarily reflect those of the European Union or the European Research Council. Neither the European Union nor the granting authority can be held responsible for them.

## Declaration of generative AI and AI-assisted technologies in the manuscript preparation process

During the preparation of this work the author(s) used ChatGPT (OpenAI, 2025) in order to improve the clarity and fluency of the writing. After using this tool/service, the author(s) reviewed and edited the content as needed and take(s) full responsibility for the content of the published article.